\documentclass[12pt]{article}
\usepackage{epic,eepic}
\begin{document}
\makeatletter
\headsep 10 mm
\footskip 10mm
\begin{center}
{\Large \bf THREE, FOUR AND FIVE-DIMENSIONAL FULLERENES
\footnote{This work was supported by the Sonderforschungsbereich 343
``Diskrete Strukturen in der Mathematik'' (Universit\"at Bielefeld) and
Russian fund of fundamental research (grant 96-01-00166).}}
\end{center}
\vskip 10 pt
\begin{center}
{ Michel Deza and Mikhail Shtogrin}
\end{center} 
\vskip 10 pt
\baselineskip 4.5 mm

\begin{abstract}
We explore some generalizations of fullerenes $F_v$ (simple polyhedra
with $v$ vertices and only 5- and -gonal faces) seen as $(d-1)$-dimensional
simple manifolds (preferably, spherical or polytopal) with only
5- and 6-gonal 2-faces. First, finite and planar (infinite) 3-fullerenes are
described. Three infinite families of spherical 4-fullerenes are
presented in Sections 4-6. The Construction A gives 4-polytopes by suitable
insertion of
fullerenes $F_{30}(D_{5h})$ into glued 120-cells. The Construction B gives
3-spheres by growing dodecahedra and barrels $F_{24}$ around of given
fullerene. The Construction C gives 4-fullerenes from special decoration of
given 4-fullerene, which add facets $F_{20}$, $F_{24}$, $F_{26}$ and
$F_{28}(T_d)$ only. Some 5-fullerenes are obtained, by a variation
of gluing of two regular tilings $5333$ of hyperbolic 4-space or of their 
suitable quotients.

\end{abstract}

\section{Introduction}

We define  here {\em d-fullerene} as a $(d-1)$-dimensional simple (i.e. 
$d$-valent) 
manifold (on any surface), such that any 2-face is 5- or 6-gon. We are
specially interested in the $d$-fullerenes, which
are {\em spherical}, i.e. homeomorphic to the $(d-1)$-sphere, and, moreover,
{\em polytopal}, i.e. convex. So, the dual of a $d$-fullerene is a
$(d-1)$-dimensional simplicial manifold, such that any $(d-3)$-face is
adjacent to five or six $(d-2)$-simplexes.

We will use the following notation:
$F_v(G)$ denotes a {\em fullerene}, i.e simple polyhedron with only
5- and 6-gonal faces, having $v$ vertices and the group of symmetry $G$.
Those polyhedra are important in Organic Chemistry; see, for example, a
monograph \cite{FM}. In 
particular, the regular dodecahedron $F_{20}(I_h)$ and the
``hexagonal barrel'' (unique $F_{24}$)
will be also denoted by $Do$ and $B_6$, respectively.

All {\em regular} (i.e. such that the group of automophisms is transitive
on $i$-faces for any $i$) partitions, which are $d$-fullerenes, are
(in classical notation of \cite {Cox}):

for $d=3$: 53 (the smallest fullerene $F_{20}$) and 63 (the
tiling of Euclidean plane by regular hexagons; it can be seen as the 
fullerene $F_{\infty}$);

for $d=4$: 533 (the 120-cell) and 633;

for $d=5$: 5333 (the tiling of hyperbolic 4-space by 120-cell).

All 4- and 5-fullerenes below are constructed from five above regular
partitions. But, perhaps, some of those constructions can be applied
also to regular star-polytopes $\frac{5}{2} 3$, the great stellated 
dodecahedron, and $\frac{5}{2} 33$; they cover 2-, 3-sphere, respectively,
7, 191 times and isomorphic, i.e. topologically equivalent, to the 
dodecahedron and the 120-cell, respectively. All other, i.e. not simple,
regular partitions with only 5- or 6-gonal 2-faces, are: hyperbolic
5m, 6m (for any $m \ge 4$); 534, 5334; 535, 5335 and 634, 635; 536. The angle
between faces of dodecahedral facet is $90^0$, $72^0$, $60^0$ for 534, 535,
536, respectively.

The purpose of this note is to expose a zoo of examples of $d$-fullerenes 
for $d=3,4,5$; we not considered embeddability of those $(d-1)$-manifolds.

\section{3-fullerenes}

Consider first 3-fullerenes. It turns out, that all {\em finite}
3-fullerenes are:

(i) {\em spherical fullerenes}, i.e.usual fullerene polyhedra on the sphere
$S^2$;

(ii) {\em toroidal} and {\em Klein bottal polyhexes}, i.e. 3-valent
partitions of the torus $T^2$ or the Klein bottle $K^2$ into hexagons;

(iii) {\em elliptic fullerenes}, i.e. 3-valent partitions of the real
projective (i.e. elliptic) plane $P^2$ into six 5-gons and some $p_6$ hexagons.

In fact,the Euler characteristic $v - e + p_5 + p_6$ of any finite
3-fullerene is $\frac{p_5}{6} \ge 0$, because $3v = 2e = 5p_5 + 6p_6$.
On the other hand, it is $2 - 2g$ for any orientable and $2-g$ for any
non-orientable surface of genus $g$. So, the only 
possible surfaces are the $S^2,T^2,P^2,K^2$ with $g = 0,1,1,2$ and
$p_5 = 12,0,6,0$, respectively; 
$v = \frac{5p_5}{3} + 2p_6$. Toroidal and Klein bottal polyhexes are
classified, for example, in \cite{Ne1}; the smallest (polyhedral) ones are
the Heawood graph (dual $K_7$) and dual $K_{3,3,3}$, respectively. The elliptic
fullerenes are exactly {\em folded} centrally-symmetric spherical fullerenes,
i.e. their $i$-faces ($i=0,1,2$) are pairs of  its antipodal $i$-faces.
The smallest one is the Petersen graph (dual $K_6$) on $P^2$, i. e. folded
Do. Actually, the centrally-symmetric fullerenes are exactly those of symmetry
$C_i, C_{2h}, D_{2h}, D_{6h}, D_{3d}, D_{5d}, T_h, I_h$; clearly, 
4 divides $v$. So, there are 30 such fullerenes with $v \le 60$:
the $F_{20}$, the $F_{32}(D_{3d})$ (one of two $F_{32}$ with 6-ring of hexagons
separating two identical blocks of six pentagons), the $F_{36}(D_{6h})$
(elongated $F_{24}$) and 3,3,3,3,6,9 of them with 40,44,48,52,56,60 vertices,
respectively. There are 107 such fullerenes with isolated pentagons and
$v \le 140$.

All compact $d$-fullerenes are finite. As an example of infinite 3-fullerene,
let us introduce {\em planar fullerenes}, i.e. 3-valent partitions of the
plane into (combinatorial) 6-gons and $p_5$ 5-gons. It turns out, that such
partitions have 
$p_5 \le 6$. For $p_5=0,1$ such 3-fullerene is unique, for any
$2 \le p_5 \le 6$ there is an infinity of them.

In fact, above criterion follows from a theorem of Alexandrov (see, for
example page 92 of \cite{Ale}), that any
complete (i.e. without boundary) metric of non-negative curvature on the
plane, can be realized as a metric of convex surface in 3-space. Consider
planar metric, with respect of which all faces of our planar fullerene
became regular polygons. Its curvature will be zero on interior points (of 
faces and edges) and non-negative on vertices. But the convex surface will
be at most a half of the sphere and so we are done. Apropos, the same argument
implies that a partition of the plane into hexagons and $p_4$ quadrangles
(respectively, $p_3$ triangles) has $p_4 \le 3$ (respectively,
$p_3 \le 2$). 

\section{4-fullerenes}

First, by analogy with above planar fullerenes, consider {\em space}
fullerenes,
i.e. 4-valent tilings of 3-space by some usual fullerenes. Examples of those
infinite 4-fullerenes are tilings used (in \cite{We} on pp. 74, 136-139,
659-664) for description of clathrate structures of some ice-like or
silicate compounds. In particular, they are space-fillings: a) by $Do$ and
$B_6$ in proportion 1:3; b) by slightly distorted $Do$ and the $F_{28}(T_d)$ in
proportion 2:1; c) by $Do$, $B_6$ and the $F_{26}$. Those and other space
fullerenes, whose cells have isolated hexagons occur also as crystaline
t.c.p. (tetrahedrally close-packing phases) of metallic alloys.  Following
simple
construction gives new space fullerene, having, as tiles, $Do$, $B_6$ and the 
$F_{36}(D_{6d})$ (elongated $B_6$). Take infinite 4-fullerenes, given in
the end of Section 5 below, and glue, by hexagons, an infinite pile of those
parallel regions, in order to fill the 3-space. Elongated $B_6$ will come
from pairs of glued $B_6$, as in Construction A, the $F_{30}(D_{5h})$
(elongated $Do$) came from $Do$.

For any finite 4-fullerene, denote by $(v, e, p, q)$ its $f$-vector,
giving the number of its $i$-faces for $i = 0,1,2,3$. The Euler characteristic,
i.e. the number $v - e + p - q$ is 0 for any finite closed 3-manifold.
Also a i-face of simple 3-manifold is the intersection of exactly 4-i 
(i+1)-faces; so, $e=2v$, $p=v+q$ and the number $p_5$ of pentagonal 2-faces
is $6q$.

All 3-faces of 
$d$-fullerenes are 3-fullerenes; we want, moreover, those fullerenes to be
close relatives of $Do$. Besides of $B_6$
(the unique next to $Do$, by the number of vertices),
two other fullerenes with isolated hexagons, the unique $F_{26}$ and unique
$F_{28}(T_d)$, are  also candidates. The duals of those
three polyhedra are known in Chemistry (under the name Frank-Kasper 
polyhedra $Z_{14},Z_{15},Z_{16}$) and in Physics, where they appear as
dislinations (rotational
defects) with respect to the vertex figure of the local icosahedral order. 
The fullerene $F_{30}(D_{5h})$ (elongated $Do$) will also appear
below, in the Construction A, as well as the $F_{32}(D_{3h})$, in
Construction D; the $F_{36}(D_{6d})$ appeared above in this Section.

Some relevant facts and analogies for 4-fullerenes are:

(i) It is well-known (see, for example, \cite{Bok}) that the boundary of the
120-cell is the unique simple equifacetted 3-sphere with (combinatorial) facet
$F_{20}$. But \cite{Sh} has shown that every 4-polytope can be approximated
arbitrary closely (in the Hausdorff distance) by a polytope whose facets
are projective images of the dodecahedron. Remind also that the 120-cell is the
universal polytope in the sense that any regular
$\le 4$-polytope, including star-polytopes, can be inscribed (vertices into
vertices) into it ( \cite{Cox}, page 269).
In Chapter 10 of \cite{FT}, it is conjectured
that 120-cell is isoperimetrically {\em best} (i.e. it has the least volume
among the 4-polytopes of unit in-radius, having 120 cells) and it is proved
that it is locally the best. See also \cite{Con}, \cite{Mi} for some
operations on the 120-cell.

(ii) Pasini \cite{Pa} proved non-existence of {\em 4-dimensional football},
i.e. equifacetted 4-fullerene with (combinatorial) facet $F_{60}(I_h)$.
Clearly, any equifacetted spherical 4-fullerene with $q$ (combinatorial)
facets $F_n$, has $v=\frac{qn}{4}$ vertices. Perhaps, 120-cell is
unique such 4-fullerene.

(iii) There exists (non-simple, of course) a tiling of Euclidean 3-space by
(116 types of polyhedra isomorphic to) the $Do$; the question about tiling of
4-space by 4-polytopes isomorphic to the 120-cell, is open (see \cite{Sc}).
$Do$, $B_6$ (actually, each {\em pentagonal i-sided prism} $B_i, i\ge 5$, see
Section 5 below) tiles alone the hyperbolic 3-space  $H^3$ with vertex
valency 6; also 535 is 12-valent tiling of $H^3$ by $Do$.

(iv) We can show, using Theorem 6 from \cite{DS}, that the skeleton of the dual
of any 4-fullerene does not embed isometrically (up to a scale) in any cubic
lattice; the 120-cell also does not embed ( \cite{DG}) this way.

The Table below presents three families of spherical 4-fullerenes, constructed
in this note.

The columns 2, 3 of the Table give
the number of vertices and 2-faces; the number of edges is $2V$, clearly. The
next five columns give the number of corresponding fullerenes among cells of
a 4-fullerene; here $F_{20}$, $F_{24}$ are $Do$, $B_6$ while $F_{28}$,
$F_{30}$ are the unique fullerenes with symmetry
$T_d$, $D_{5h}$, respectively, and given number of vertices.
The last column $F'$ gives: the number of 3-cells $F$, when
it is a fullerene in the Construction B, and the number of 3-cells {\em in}
$F$, when it is a 4-fullerene in the construction C. The symbols $v$,
$p=(p_5,p_6)$, $q$ denote the number of vertices, of 2-faces and (for the
Construction C) of cells of $F$. In particular, $C_1(120-cell)$ has, as cells,
$2p_5 + q =1560$ dodecahedra and $v = 600$ fullerenes $F_{28}$.

\vspace{0.3cm}
{\bf Table. 
 f-vectors of some finite 4-dimensional fullerenes}
\[\begin{array}{|c||c||c||c|c|c|c|c||} \hline
 &V&P & F_{20} &F_{24} &F_{28}&F_{30}&F'\\ \hline \hline
120-&600&720&120&0&0&0& \\ \hline  
A_i&560i+40&666i+54&94i+26&0&0&12i-12& \\
 \hline  
B(F)&30v&\frac{71v}{2}+10&\frac{7v}{2} +48&2v - 40&0&0&2 \\ \hline        
C_1(F)&20v&20v+3p&2p_5&2p_6&v&0&q    \\ \hline \hline
\end{array} \]  

\section{Construction of polytopal 4-fullerenes $A_i$}

\medskip
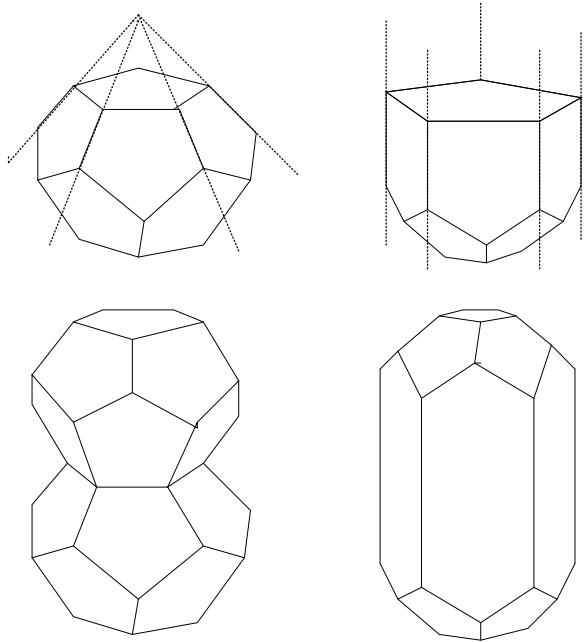
\begin{center}
\setlength{\unitlength}{0.00041667in}
\begingroup\makeatletter\ifx\SetFigFont\undefined%
\gdef\SetFigFont#1#2#3#4#5{%
  \reset@font\fontsize{#1}{#2pt}%
  \fontfamily{#3}\fontseries{#4}\fontshape{#5}%
  \selectfont}%
\fi\endgroup%
{\renewcommand{\dashlinestretch}{30}
\begin{picture}(7299,8139)(0,-10)
\path(312,3387)(837,4062)(1587,3837)(1587,3162)
\path(2937,3312)(2487,4062)(1587,3837)
\path(837,4062)(1212,4212)(2112,4212)(2487,4062)
\path(2487,2262)(3087,1662)(3012,1062)(2487,1212)
\path(762,2262)(312,1737)(312,1062)(837,1212)
\path(312,1062)(837,312)(1587,87)(1662,537)
\path(1587,87)(2412,237)(3012,1062)
\path(4737,3462)(4737,987)
\path(5262,3087)(5262,687)
\path(6687,3087)(6687,687)
\path(7212,3462)(7212,987)
\path(5262,3087)(5937,3537)(6687,3087)
\path(5262,687)(6012,237)(6687,687)
\path(4737,987)(4962,537)(5262,687)
\path(7212,987)(6987,537)(6687,687)
\path(4962,537)(5487,87)(6012,12)
	(6612,162)(6987,537)
\path(6012,237)(6012,12)
\path(4737,3462)(4962,3687)(5262,3087)
\path(7212,3462)(6912,3762)(6687,3087)
\path(4962,3687)(5487,4137)(6012,4062)
	(5937,3537)(6012,3537)
\path(6912,3762)(6462,4137)(6012,4062)
\path(5487,4137)(5787,4212)(6237,4212)(6462,4137)
\path(837,7062)(387,6537)(387,5862)(912,6012)
\path(387,5862)(912,5112)(1662,4887)(1737,5337)
\path(762,2262)(312,3012)(312,3387)(837,2787)
\path(1662,4887)(2487,5037)(3087,5862)
\path(762,2262)(1137,1962)(2037,1962)
	(2487,2262)(2937,2862)(2937,3312)
	(2412,2787)(2037,1962)(2487,1212)
	(1662,537)(837,1212)(1137,1962)
	(837,2787)(1587,3162)(2412,2712)(2412,2787)
\path(5337,5487)(6087,5037)(6762,5487)
\path(4812,5787)(5037,5337)(5337,5487)
\path(7287,5787)(7062,5337)(6762,5487)
\path(6087,5037)(6087,4812)
\path(837,7062)(1212,6762)(2112,6762)
	(2562,7062)(1662,7287)(837,7062)
\path(6012,7137)(4812,6987)(5337,6612)
	(6762,6612)(7287,6912)(6012,7137)
\path(4812,6987)(4812,5787)
\path(5337,6612)(5337,5487)
\path(6762,5487)(6762,6612)
\path(7287,6912)(7287,5787)
\dottedline{45}(1587,7887)(1662,7887)(537,5037)
\dottedline{45}(1662,7962)(3687,5937)
\dottedline{45}(4812,7887)(4812,5037)
\dottedline{45}(5337,7512)(5337,4737)
\dottedline{45}(6762,7512)(6762,4737)
\dottedline{45}(7287,7737)(7287,5112)
\dottedline{45}(6012,7137)(6012,8112)
\path(5037,5337)(5562,4887)(6087,4812)
	(6537,4962)(7062,5337)
\path(5337,6612)(6762,6612)(7287,6912)
	(6012,7137)(4812,6987)(5337,6612)
\dottedline{45}(1662,7962)(2937,4962)
\dottedline{45}(1662,7962)(12,6087)(12,6162)
\path(1212,6762)(912,6012)(1737,5337)
	(2487,6012)(2187,6762)(2112,6762)
\path(2562,7062)(3162,6462)(3087,5862)(2487,6012)
\end{picture}
}

\end{center}
\medskip

It will be a 4-dimensional analogue of the following simple construction of
the $i$-layered dodecahedron $F_{20+10i}$; see Figure above for such $ F_{30}$.
Stellate a face $t$ of $Do$ (i.e. extend face-planes of its 5 neighbors until
their intersection; so we got a 5-pyramid on the face). Then do a projective
transformation, sending the apex of 5-pyramid to infinity so that the
5-pyramid became right regular 5-prism. The image of our regular dodecahedron
will be inscribed in the continuation of above 5-prism. The face $t$
became larger and its opposite became smaller, but they both remain regular
5-gons; all other faces became irregular. Only one of the six 5-axes of
symmetry of $Do$ will remain. Take the mirror reflection of
such modified dodecahedron on the face-plane of $t$. Two such dodecahedra
glued by the ``large'' regular face, obtained from $t$, form the convex
3-polytope $F_{30}(D_{5h})$. It has exactly two regular 5-gonal faces:
``small'' ones from both modified dodecahedra. On each of them we can continue
the same procedure and get general $i$-layered dodecahedron $F_{20+10i}$ with
symmetry $D_{5d}$ for even $i > 0$ and $D_{5h}$ for odd $i$. (This tube is the
dual of 2-capped pile of $i$ 5-anti-prisms.)

Apply same procedure to the 120-cell in 4-space. Stellate a dodecahedral
face $t$ until we get a pyramid on it. By a projective transformation, sending
the apex to infinity, it will be transformed into a right prism, having $Do$ as
a base. The 120-cell will be modified: $t$ became larger, its opposite became
smaller, but both remain regular; all other dodecahedral cells became
irregular. Take the mirror reflection of the modified 120-cells on the 3-space,
containing $t$; we get (from two modified 120-cells, glued by the ``large''
regular dodecahedron) the convex 4-polytope $A_1(120-cell)$. It has exactly
two faces $Do$, ``small'' $Do$'s of two modified 120-cells, the other
dodecahedra
are irregular. The continuation of this procedure on each of ``small'' $Do$
gives the 4-polytope $A_i(120-cell)$. See its $f$-vector in the Table above;
exactly $30i-30$ its 2-faces are hexagons.
 
We can apply the construction A to any {\em non-exposed} dodecahedral cell of a
$A_i(120-cell)$, i.e. having only dodecahedral neighbors: we obtain 3-spheres,
but now there is no guarantee of convexity. When operation A is applied to
several non-exposed dodecahedra, no two of them should have the same
dodecahedral neighbor. In order to enumerate such possibilities,
the solution of the following extremal problem will be of interest
(we give it in dual form for the 600-cell):
estimate the maximal number of vertices in the 600-cell with all pairwise
distances (in the skeleton, having diameter 5) at least 3. It is at most
9, clearly, and  at least 6: take 3 suitable vertices on each of two 10-gons
(among all 12), which lie in two orthogonal planes. 

In Chapter 4 (Sections 1.7 and 1.8) of \cite{SM}, four 4-fullerenes,
having each 144 dodecahedral cells  and $12k$ cells $B_6$ for $k=2,3,4,6$, are
constructed from the 600-cell, by inverting the Hopf map of the 3-sphere on
the 2-sphere. Also in Chapter 7 (Sections 2.7, 3 and 4.2) of \cite{SM} crystal
agregats are given, which can be used to construct 4-fullerenes.

\section{Construction B of spherical 4-fullerenes}

Fix a 3-fullerene $F$ with $v$ vertices, $p = \frac{v}{2} + 2$ faces and
$e = \frac{3v}{2}$ edges.
From an interior point $o$ of $F$ take on the ray through each vertex $b$
a point $b'$ with distance $d(o,b') = d(o,b) + 1$. Put on each face of $F$
dodecahedra $Do$ on 5-gons and barrels $B_6$ on 6-gons, so that their lateral
sides coincide. (Always in this construction $Do$ and $B_6$ are
combinatorial.) 
We got {\em 1-corona}: $F$ itself and $p$ polyhedra of the {\em 1-st floor}.
The 
surface of 1-corona consists of $p$ {\em 1-anti-faces}, i.e. opposite ones to
the faces of $F$ and others, which are organized in $v$ 3-hedral triples of
5-gons with the central vertex $b'$ (for each of $v$ vertices $b$ of $F$).
Put $v$
new $Do$ into those $v$ 3-hedral angles, one $Do$ for each. We got a 2-corona
with the 2-nd floor, consisting of $v$ dodecahedra. Each of them is adjacent to
1-corona in 3 faces (of its 3-hedral angle) and to 3 neighbors on the 2-nd
floor; so that the6 remaining faces are free. Each of the 12 5-gonal (or
$p - 12$
hexagonal) 1-anti-faces is incident to 5 (or 6,respectively) dodecahedra of
the 2-nd floor. Those 5 (or 6) 5-gons form a half-dodecahedron (or a
half-barrel, respectively). Add for   
each of them the second half in order to obtain $p$ new polyhedra; they form
3-rd floor. We got 3-corona. The surface of 3-corona consists of $p$
2-anti-faces (i.e. the faces, opposite to 1-anti-faces) and $e$
{\em quadruples}, i.e. two edge-adjacent 5-gons and two other 5-gons,
edge-adjacent to the first two via each vertex of the edge of their adjacency.
First two 5-gons are from the surface of 2-corona, two others are from the
surface of the 3-corona. Take now two copies of 3-corona. (Remind, that each
i-corona is a 3-ball in 3-space.) Now
we will join them in 4-space, putting between them $e$ new dodecahedra, which
will form the 4-th floor for each copy. Also the corresponding 2-anti-faces
of them will coincide. Each $Do$ of the 4-th floor is incident to each copy of
3-corona by a quadruple and to four neighbors on the 4-th floor.   

Clearly, $B(Do)$ is the 120-cell itself and $B(B_6)$ consists only of two
(combinatorial types of) fullerenes $Do$ and $B_6$.

In fact, the construction B can be similarly applied to
 any simple 3-polytope with, say, $v$ vertices and any given
{\em p-vector} $(p_3, p_4,...)$, where $p_i$ is the number of $i$-gonal faces
for any $i \ge 3$. The above construction will give simple 3-sphere with $30v$
vertices, $60v$ edges, $\frac{71v}{2} + 10$ 2-faces (including $5p_i$
$i$-faces for each $i$, except 5) and $\frac{11v}{2} + 10$
cells, including two original 3-polytopes, $4p_i$ $i$-gonal barrels $B_i$
($p_i$ on both 1-st and both 3-rd floors) and $\frac{7v}{2}$
dodecahedra ($v$ on each of both 2-nd floors and $\frac{3v}{2}$, i.e. the
number $e$ of edges, on the common 4-th floor). Note that $B_5$ is $Do$ and so,
$4p_5 + \frac{7v}{2}$ is the number of all dodecahedra.

It looks difficult to decide when the construction B leads to a 4-polytope
(i.e. a convex 3-sphere) even when applied to such polyhedra as the regular
tetrahedron, the cube or the barrel $B_i$. 
$B_3$ is a cube with two opposite vertices truncated; $B_4$ is the dual of
2-capped 4-anti-prism (one of all 8 polyhedra whose faces are regular
triangles). $B_3$ is called in  \cite{BS} {\em D\"urer octahedron};
among constructions, given in \cite{BS}, there are two simple equifacetted
3-manifolds: one with 10 facets $B_3$ and
one with 26 facets $B_4$ (the first one is a non-polytopal 3-sphere). 

Construction B can be applied also to any simple partition of Euclidean or
hyperbolic plane. In particular, taking, as original $F$, the ``graphite''
$63=F_{\infty}$ (Euclidean plane, partitioned by regular hexagons), one get
a region of 3-space, bounded by two parallel planes and filled  by several
layers of $Do$ and $B_6$. 

\section {Construction C of 4-fullerenes $C_j(F)$}

We give this construction in dual terms of general simplicial 3-manifold;
it was inspired by \cite{SM1}.
Apply following four operations to the simplicial 3-manifold $F^*$, which is
the dual to given simple 3-manifold $F$:

1) Transform each edge into a 4-path of three edges, by adding two new
  vertices on each edge, and subdivide each tetrahedron, using the new edges, into
  four tetrahedra and one truncated tetrahedron.

2) By projecting all faces from an interior point, subdivide each 
  truncated tetrahedron into four tetrahedra and four 6-pyramids.

3) Glue each two 6-pyramids having a common base into a 6-bipyramid (cf. two
  4-pyramids, glued into the octahedron in transition to the f.c.c.
lattice $A_3)$.

4) Subdivide each 6-bipyramid into six tetrahedra with common edge, linking
  its apexes.

Denote the obtained simplicial complex by $C_1(F)$; iterating the above
procedure $j$
times produces $C_j(F)$. If $F$ has $v$ vertices, $p$ 2-faces (including $p_5$
5-gonal and $p_6$ 6-gonal ones), $q$  cells , then
$C_1(F)$ has $20v+3p$ 2-faces (including $2v+3p_6$ hexagons) and only following
cells: all cells of $F$ plus $2p_5$ dodecahedra,
$2p_6$ hexagonal barrels $B_6$ and $v$ fullerenes $F_{28}(T_d)$.    
So, if $F$ is a 4-fullerene (for example, the 120-cell or one obtained by the
above constructions A or B), then any $C_j(F)$ is also 4-fullerene. 

If the original simplicial manifold $F^*$ is spherical, then its $j$-th simplicial subdivision, described above, is also spherical. But the question of
preserving the convexity is difficult. Operations 1), 2), 4) could be arranged
in order to preserve it. (For example, chosen interior points of the
tetrahedra should be moved ``out'' within 4-th dimension in order to get edges
between neighbors, then suitable two points {\em around} of each edge should
be found and so on.) But the operation 3) can destroy the convexity.
Moreover, four above topological operations can be seen separately, which is
not the case of their metrical contreparts.

The dualization of another decoration of 600-cell, given in \cite{Mo} and 
Section 3.4 of \cite{SM1}, produces from it another infinite family $D_j, 
j \ge 0,$ of
spherical 4-fullerenes, having $61^j \times 600$ vertices and, as cells, only
dodecahedra, $B_6$ and the fullerene $F_{32}(D_{3h})$.
Similarly to the construction $C_j(F)$, one can generalize it to a construction
(say, $D_j(F)$), which gives an infinite
family of 4-fullerenes (having, besides of cells of $F$, only cells $B_6$ and
$F_{32}(D_{3h})$) starting from a 4-fullerene $F$. A mixed construction
(choosing suitably operation C or D on
each step) gives asymptotically non-periodic 4-fullerenes, having, besides of
cells of $F$, only cells $B_6$, $F_{28}(T_d)$ and $F_{32}(D_{3h})$.

Both decorations of 600-cell, leading to $C_j(120-cell)$ and $D_j(120-cell)$,
are given in \cite{SM1} as two examples of {\em inflation} method, which can
be used
in order to get other spherical $d$-fullerenes. Roughly, it consists of finding
out, in a simplicial $d$-polytope $F^*$, a suitable ``large'' $(d-1)$-simplex,
containing an integer number of ``small'' (fundamental) simplices; this number
was 20, 61 for Constructions $C_1(F)$, $D_1(F)$ with $F$ being $120-cell$. The
decoration of $F^*$ comes by ``barycentric homothety'' (suitable projection
of the ``large''simplex on the ``small''one) as the orbit of new points under
the symmetry group. 

\section{Quotient d-fullerenes and polyhexes} 

The terms of this Section are, for example, from \cite{Thu}, \cite{Rag}.

The Poincarè {\em dodecahedral space} (\cite{Poi})
is, actually, a 4-fullerene with $f$-vector
$(v,e,p,c)=(5,10,6=p_5,1)$. It comes from $Do$ by gluing of its opposite faces
with $ \frac{1}{10} $ right-handed rotation. The Poincarè fullerene is
{\em not} spherical, but it is {\em locally} spherical (i.e. of constant
positive curvature); it has same Betti
groups as $S^3$, but its fundamental group is not trivial (it has order 120).
Actually, the Poincarè fullerene is the {\em quotient} 3-manifold of $S^3$ by
the binary icosahedral group $I_h$ of order 120; so its $f$-vector is just
$\frac{1}{120}$ of that one of 120-cell.

Following three  compact quotients 3-manifolds have only 5-gonal and (only
for the third one, called Löbell space) 6-gonal 2-faces, but they are not
4-fullerenes: (i) Seifert-Weber hyperbolic space with $f$-vector 
$(1,6,6=p_5,1)$, which is the quotient of $535$ by the torsion-free group,
having $Do$ as the fundamental domain; (ii) the quotient of $534$
(or of 6-valent tiling of $H^3$ by $B_6$) by a torsion-free group,
having eight right-angled $Do$ (or $B_6$) as the fundamental domain, with
$f$-vector $(20,60,48=p_5,8)$ (or $(24,72, 48+8=p_5+p_6,8)$).

In general, a theorem of Selberg (see Lemma 8 in \cite{Sel}, \cite{Bor}, 
6.11-6.16 in \cite{Rag}) gives,
that if a discrete group of motions of a symmetric space has compact
fundamental domain, then it has a torsion-free normal subgroup of finite
index. So, a quotient of a $d$-fullerene by a such group of symmetry, is a
finite $d$-fullerene.

This method gives, for example, compact 4-fullerenes as quotients of the
space fullerenes, given in the beginning of Section 3, and, on 
cylinder or half-cylinder, of 4-fullerenes, given in the end of Section 5.

Consider now {\em polyhexes}, i.e. $d$-fullerenes, having only 6-gonal
2-faces. The regular tilings $63=F_{\infty}$ and $633$ are examples of such
infinite simply connected 3- and 4-fullerenes. All quotient surfaces of the
Euclidean plane $R^2$ by discontinuous and fixed point free group of
isometries, are: torus $T^2$, cylinder, its twist (Möbius surface) and Klein
bottle $K^2$; the group is generated, respectively, by two translations, a
translation, a glide reflection and by a translation {\em and} glide
reflection. So, those four types of quotients of $63$ are polyhexes; they are 
not simply connected.

Smallest quotient polyhexes on $T^2$ realize graphs
$C_6$, $K_{3,3}$ (two embeddings), 3-cube with $p_6=1,3,4$, respectively. But
such polyhexes, the set of all faces, edges and vertices of which form a
topological complex (i.e. no loops and double edges in the graph and the
intersection of any two faces is an edge, a vertex or empty), exist if and only
if $p_6 \ge 7$. On the other hand, the ``greatest''polyhex $633$ (the convex
hull of the vertices of $63$, realized on a {\em horosphere}, i.e. on a sphere
with center at infinity) is an example of not
{\em cocompact} (i.e. with not compact fundamental domain), but {\em cofinite}
(i.e. of finite volume) 4-fullerene.

\section{5-fullerenes}

The regular tiling $5333$ of hyperbolic 4-space by 120-cells
is an infinite 5-fullerene: all its 2-faces are 5-gons.

The following is an infinite family of 5-fullerenes, having both
5-gonal and 6-gonal 2-faces.
Take two copies of the tilings $5333$ and glue them in some pairs of
corresponding 120-cells. Delete now from the manifold the interiors of
those 120-cells. For each of them, any corresponding pair (from both
$5333$) of neighboring 120-cells glue in a 4-polytope $A_1$, described in
the Section 4. If the tilings are glued in only one 120-cell, the
4-manifold is the direct product of the 3-sphere and the Euclidean line;
so it is simply connected.

Now, using Selberg theorem, as in previous Section, we get a {\em finite}
5-fullerene as quotient of $5333$ by its symmetry group. It is a compact
4-manifold, partitioned into a finite number of 120-cell's. The same gluing
as above, can be applied to it, in order to produce new 5-fullerenes.


{\it 
LIENS-DMI, Ecole Normale Sup\'erieure, 45 rue d'Ulm, 75005 Paris, FRANCE.
e-mail: deza@dmi.ens.fr

Steklov Mathematical Institute, Gubkin str. 8, 117966 Moscow GSP-1, RUSSIA.
e-mail: stogrin@mi.ras.ru }

\begin{thebibliography}{99}

\bibitem[Ale58]{Ale}
A.D.Alexandrov, {\em Convexe Polyheder},
Akademie-Verlag, Berlin, 1958.

\bibitem[BoS95]{BS}
J. Bokowski and P. Schubert,
{\em Equifacetted 3-spheres as Topes of Nonpolytopal Matroid Polytopes},
Discrete Comput. Geometry {\bf 13} (1995) 347--361.

\bibitem[Bok95]{Bok}
J.Bokowski, {\em On the construction of equifacetted 3-spheres},
in N.L.White (ed.), Invariant Methods in Discrete and Comput. Geometry,
Kluwer Acad. Publ. (1995) 301--312.

\bibitem[Bor63]{Bor}
A.Borel, {\em Compact Clifford-Klein forms of symmetric spaces},
Topology {\bf 2} (1963) 111--122.

\bibitem [Con67]{Con}
J.H.Conway, {\em Four-dimensional Archimedean polytopes}, in
Proc. Colloquium on Convexity, Copenhagen 1965, Kobenhavns Univ. Mat. Institut
(1967) 38--39.

\bibitem[Cox73]{Cox} H.S.M.Coxeter,
{\em Regular Polytopes},
 3rd ed., Dover, New York, 1973.
  
\bibitem[DGr97]{DG}
M.Deza and V.P.Grishukhin, {\em The skeleton of the 120-cell is not 5-gonal},
Discrete Math. {\bf 165/166} (1997) 205--210.

\bibitem[DSt97]{DS}
M.Deza and M.I.Shtogrin, {\em Embedding of skeletons of Voronoi and Delone
partitions into cubic lattices},
in Voronoi's impact on modern science, Book 2, Institute of Mathematics, Kyiv
(1998) 80--84.

\bibitem[FoMa95]{FM}
P.W.Fowler and D.E.Manolopoulos, {\em An Atlas of Fullerenes},
Clarendon Press, Oxford, 1995.

\bibitem[Miy90]{Mi}
K.Miyazaki, {\em Primary Hyper-geodesic Polytopes},
Int. J. of Space Structures {\bf 5}(3-4) (1990) 309--323.

\bibitem[Mo83]{Mo}
R.Mosseri, {\em These d'Etat}, Orsay, 1983.

\bibitem[Neg85]{Ne1}
S.Negami, {\em Uniqueness anf faithfulness of embedding of graphs into
surfaces},
Doctor thesis, Tokyo Institute of Technology, 1985.

\bibitem[Pas98]{Pa}
A.Pasini, manuscript, January 1998.

\bibitem[Poi04]{Poi}
H.Poincaré, {\em Cinquième complèment à l'analysis situs},
Rend. Circ. Mat. Palermo {\bf 18} (1904) 45--110, or {\em Oeuvres}, t.6
435--498.

\bibitem[Rag72]{Rag}
M.S.Raghunathan, {\em Discrete Subgroups of Lie Groups},
Springer-Verlag, Berlin, 1972. 

\bibitem[SaM85]{SM1}
J.F.Sadoc and R.Mosseri, {\em Hierarchical interlaced networks of
disclinations lines in non-periodic structures},
J. Physique {\bf 46} (1985) 1809--1826.

\bibitem[SaM97]{SM}
J.F.Sadoc and R.Mosseri, {\em Frustration géométrique},
Eyrolles, Paris, 1997.

\bibitem[Sel60]{Sel}
A.Selberg, {\em On discontinuous groups in higher-dimensional symmetric
spaces}, in {\em Contributions to Function Theory}, TIFR, Bombay (1960)
147--164. 

\bibitem[She66]{Sh}
G.C.Shephard, {\em Approximation by polytopes with projectively regular
faces}, Mathematika, London {\bf 13} (1966) 189--195.

\bibitem[Sch84]{Sc}
E.Schulte, {\em Non-tiles and non-facets for the Euclidean space, spherical
complexes and convex polytopes},
J. Reine Angw. Math. {\bf 352} (1984) 161--183.

\bibitem[Thu97]{Thu}
W.P.Thurston, {\em Three-Dimensional Geometry and Topology}, Princeton Univ.
Press, Princeton, New Jersey,1997.

\bibitem[FTo64]{FT}
L.Fejes T\'oth, {\em Regular Figures},
Pergamon Press, Oxford, 1964.

\bibitem[Wel84]{We}
A.F.Wells, {\em Structural Inorganic Chemistry},
5-th ed., Oxford, 1984.

\end{thebibliography}
\end{document}